\documentclass{article}

\usepackage[english]{babel}
\usepackage{amsmath}
\usepackage{amsthm}
\usepackage{amssymb}

\newtheorem{thm}{Theorem}

\newtheorem{lem}{Lemma}
\newtheorem{cor}{Corollary}

\theoremstyle{definition}
\newtheorem{Def}{Definition}

\theoremstyle{remark}

\theoremstyle{remark}

\def\cE{\mathcal E}

\def\cM{\mathcal M}

\def\cP{\mathcal P}

\def\cH{\mathcal H}

\def\RR{\mathbf R}
\def\bR{\mathbf R}

\def\NN{\mathbf N}

\def\Im {\mbox{\rm Im}\,}
\def\Ker{\mbox{\rm Ker}\,}
\def\la{\lambda}

\begin{document}

\title{On polynomial solutions of generalized Moisil-Th\'eodoresco systems and\\ Hodge-de Rham  systems}

\author{Richard Delanghe, Roman L\' avi\v cka and Vladim\'ir Sou\v cek}

\date{}

\maketitle

\begin{abstract}
The aim of the paper is to study
relations between polynomial solutions of generalized Moisil-Th\'eodoresco (GMT) systems and polynomial solutions of
Hodge-de Rham systems and, using these relations, to describe polynomial solutions of GMT systems.
We decompose the space of homogeneous solutions of GMT system of a given homogeneity into irreducible pieces
under the action of the group $O(m)$ and we characterize individual pieces by their highest weights and
we compute their dimensions.
\end{abstract}


%%%%%%%%%%%%%%%%%%%%%%%%%%%%%%%%%%%%%%%%%%%%%%%%%%%%%%%%%%%%%%%%%%%%%%%%%%%%%%%%%%%
%%%%%%%%%%%%%%%%%%%%%%%%%%%%%%%%%%%%%%%%%%%%%%%%%%%%%%%%%%%%%%%%%%%%%%%%%%%%%%%%%%%%

\section{Introduction}

In its classical form Clifford analysis is the study of properties of solutions of the Dirac operator $D$
acting on functions defined on $\bR^m$ with values in the corresponding Clifford algebra $\bR_{0,m}.$  It is, however, often important (and interesting) to consider  special types of solutions obtained by considering functions taking values in
suitable subspaces of $\bR_{0,m}.$ To describe some of such important cases, we shall use the language of differential forms.

Let $G$ be an open subset of the Euclidean space $\bR^m$ and let $$\cE(G)=\bigoplus_{s=0}^m\cE^s(G)$$ be the space of (smooth) differential forms
on $G.$
It is well known that the Clifford algebra $\bR_{0,m}$ is isomorphic (as a vector space)
with the Grassmann algebra $\Lambda^*(\bR^m).$  Hence the space of (smooth) Clifford algebra valued functions on $G$ can be identified with the space $\cE(G)$ of (smooth) differential forms  on $G.$ As was explained in detail in \cite{BDS},
the Dirac operator $D$ corresponds under this identification to the operator $d+d^*$
acting on the space of differential forms, where $d$ and $d^*$ are the standard
de Rham differential and its adjoint, respectively.

In recent years, there was a growing interest in the study and better understanding
of properties of solutions of generalized Moisil-Th\'eodoresco systems (see \cite{GM,BRDS,del07,del,
IC,BG}), which are defined as follows (see \cite{BRDS}).

\begin{Def}
Assume that  $r,p$ and $q$ are non-negative integers such that $p\leq q$ and $r+2q\leq m.$
Furthermore, denote by $\cE^{(r,p,q)}(G)$ the subspace of $\cE(G)$ determined by
$$\cE^{(r,p,q)}(G)=\bigoplus_{j=p}^q\cE^{r+2j}(G).$$
A generalized Moisil-Th\'eodoresco system of type $(r,p,q)$ (GMT system for short) is then defined as the
homogeneous system obtained by restricting the
operator $d+d^*$ to the space $\cE^{(r,p,q)}(G),$ i.e.
$$(d+d^*)\;\omega=0,\ \omega\in\cE^{(r,p,q)}(G).$$
\end{Def}

Notice that if $\omega\in\cE^{(r,p,q)}(G)$ is written out as
$$\omega=\sum_{j=p}^q\omega^{r+2j}$$ with $\omega^{r+2j}\in\cE^{r+2j}(G),$
then the equation $(d+d^*)\;\omega=0$ means that
\begin{eqnarray}\label{dirac2}
d^*\omega^{r+2p} &=& 0,\nonumber\\
d\omega^{r+2j}+d^*\omega^{r+2j+2} &=& 0,\ j=p,\ldots,q-1,\\
d\omega^{r+2q} &=& 0.\nonumber
\end{eqnarray}

Some special cases of GMT systems are well known and well understood. It is possible to choose
$(r,p,q)$ in such a way that $\cE^{(r,p,q)}(G)$ is equal to the space $\cE^+(G)$ (resp. $\cE^-(G)$)
of all even (resp.\ odd) forms. In these cases, the corresponding GMT system is (equivalent to) the classical Dirac equation for functions with values in the even (resp.\ odd) part of the Clifford algebra.
Properties of these functions were carefully studied in classical Clifford analysis.

Another very important special case is the GMT system of type $(r,0,0).$ In this case,
the space $\cE^{(r,p,q)}(G)$ reduces to the space $\cE^r(G)$ of forms of pure degree $r$ and
the corresponding GMT system coincides with the Hodge-de Rham system
$$
d\omega=0,\ d^*\omega=0,\ \ \ \omega\in\cE^r(G).
$$
Solutions of the Hodge-de Rham system will play a key role in our study of solutions of GMT systems.
As discussed in \cite{BRDS}, it is important to understand the structure of homogeneous polynomial solutions of these systems.

The aim of the paper is to study
relations between polynomial solutions of GMT systems and polynomial solutions of
Hodge-de Rham systems. Using these relations and results obtained by Y. Homma in \cite{H} on solutions of Hodge-de Rham systems, 
it is  then possible to describe fully polynomial solutions of GMT systems.
To describe them in more details, we shall introduce a suitable notation.

Denote by $\cP$ the space of real-valued polynomials in $\RR^m$ and by $\cP_k$ the space of polynomials of $\cP$ which are homogeneous of degree $k.$ Moreover, let $\Lambda^*(\RR^m)$ and $\Lambda^s(\RR^m)$ stand for the exterior algebra and the space of $s$-vectors over $\RR^m,$ respectively. Of course,
$$\cP=\bigoplus_{k=0}^{\infty}\cP_k\mbox{\ \ and\ \ }\Lambda^*(\RR^m)=\bigoplus_{s=0}^m\Lambda^s(\RR^m).$$
Let us now introduce the following spaces of differential forms with polynomial coefficients: $\cP^*_k=\cP_k\bigotimes_{\RR}\Lambda^*(\RR^m),$ $\cP^s_k=\cP_k\bigotimes_{\RR}\Lambda^s(\RR^m)$ and, finally,
$$\cP^{(r,p,q)}_k=\bigoplus_{j=p}^q\cP^{r+2j}_k.$$

In what follows, we shall mainly study  the space
\begin{equation}\label{MT}
MT^{(r,p,q)}_k=\{P\in\cP^{(r,p,q)}_k:\ (d+d^*)P=0\}
\end{equation}
of polynomial solutions of the generalized Moisil-Th\'eodoresco system of homogeneity $k.$
We shall relate them to the spaces
\begin{equation}\label{Hsk}
H^s_k=\{P\in\cP^s_k:\ dP=0,\ d^*P=0\}
\end{equation}
of polynomial solutions of the Hodge-de Rham  system of homogeneity $k$ and use the results
obtained by Y. Homma in \cite{H} for $H^s_k.$ 

In Section 2, we decompose $MT^{(r,p,q)}_k$ into a sum of pieces isomorphic to spaces $H^s_k$ of polynomial solutions  of various Hodge-de Rham systems (see Theorem \ref{decomposition}). In such a way,
we shall be able to characterize irreducible components of these spaces under the action
of the group $O(m)$ and count the dimension of the space of polynomial solutions of given homogeneity.
  In Section 3,
we shall review results from \cite{H} on the decomposition of the space $\Ker_k^s\Delta$ of harmonic polynomials of homogeneity $k$ with values in
$\Lambda^s(\RR^m)$ into four different pieces and the description of the spaces $H^s_k.$  In the last section, we shall add an alternative description of the four pieces in the decomposition of $\Ker_k^s\Delta.$

 %%%%%%%%%%%%%%%%%%%%%%%%%%%%%%%%%%%%%%%%%%%%%%%%%%%%%%%%%%%%%%%%%%%%%%%%%%%%%%%%%%%%%%%%%%%%%%%%%%%%%%

\section{The decomposition of $MT^{(r,p,q)}_k$}

In this part, we shall describe a decomposition of the spaces  $MT^{(r,p,q)}_k$
into a~direct sum of spaces isomorphic to spaces of solutions of various Hodge-de Rham systems.
It will give us (as we shall see later) a decomposition of the spaces  $MT^{(r,p,q)}_k$
into a sum of irreducible components under the $O(m)$-action.
We shall need the following refined version of the Poincar\'e Lemmas
(see also \cite{BDS} and \cite{del}).

\begin{lem}\label{poincare}
The following properties hold:\\
(i) For $P^s_k\in\cP^s_k$ with $s>0,$ $dP^s_k=0$ if and only if there exists $P^{s-1}_{k+1}\in\cP^{s-1}_{k+1}$ such that $d^*P^{s-1}_{k+1}=0$ and $dP^{s-1}_{k+1}=P^s_k.$\\
(ii) For $P^s_k\in\cP^s_k$ with $s<m,$ $d^*P^s_k=0$ if and only if there exists $P^{s+1}_{k+1}\in\cP^{s+1}_{k+1}$ such that $dP^{s+1}_{k+1}=0$ and $d^*P^{s+1}_{k+1}=P^s_k.$
\end{lem}

\noindent
Notice that $H^s_0=\cP^s_0$ and $MT^{(r,p,q)}_0=\cP^{(r,p,q)}_0.$

Now we can state the main result of this section.

\begin{thm}{\label{decomposition}} Let $k\in\NN.$
Denote by $\Phi$ the restriction of the operator $d$   to the space $MT^{(r,p,q)}_k.$
Then
$$\Ker \Phi=\bigoplus_{j=p}^qH^{r+2j}_k\mbox{\ \ and\ \ }\Im \Phi=\bigoplus_{j=p}^{q-1}H^{r+2j+1}_{k-1}.$$
Moreover, we have that
 $$MT^{(r,p,q)}_k\simeq\Ker\Phi\oplus\Im\Phi.$$
 \end{thm}

\begin{proof}
Let $P_k\in MT^{(r,p,q)}_k,$ that is,
$P_k=P^{r+2p}_k+P^{r+2p+2}_k+\cdots+P^{r+2q}_k$ with $P^{r+2j}_k\in\cP^{r+2j}_k.$ By (\ref{dirac2}) we have that
\begin{eqnarray}\label{dirac3}
d^*P^{r+2p}_k &=& 0,\nonumber\\
d P^{r+2j}_k &=& -d^*P^{r+2j+2}_k,\ j=p,\ldots,q-1,\\
d P^{r+2q}_k &=& 0.\nonumber
\end{eqnarray}
For $j=p,\ldots, q-1,$ put $P^{r+2j+1}_{k-1}=dP^{r+2j}_k.$ Obviously, $P^{r+2j+1}_{k-1}=-d^*P^{r+2j+2}_k$
and $\Phi(P_k)=P_{k-1}^{r+2p+1}+P^{r+2p+3}_{k-1}+\cdots+P^{r+2q-1}_{k-1}.$
Moreover, by virtue of (\ref{dirac3}), we have that
$$\Ker \Phi=\bigoplus_{j=p}^qH^{r+2j}_k.$$
Using the fact that $d^2=0$ and $(d^*)^2=0,$ we get that
$$\Im \Phi\subset\bigoplus_{j=p}^{q-1}H^{r+2j+1}_{k-1}.$$
To show the opposite inclusion, consider an arbitrary form $$P_{k-1}=P_{k-1}^{r+2p+1}+P^{r+2p+3}_{k-1}+\cdots+P^{r+2q-1}_{k-1}$$
with $P_{k-1}^{r+2j+1}\in H^{r+2j+1}_{k-1}.$
Then it is sufficient to find a~form $P_k\in MT^{(r,p,q)}_k$ such that $\Phi(P_k)=P_{k-1}.$
By using Lemma \ref{poincare}, we shall construct such a~$P_k$ as follows.\\
(i) For $dP_k^{r+2p+1}=0$ we can find $P_k^{r+2p}\in\cP_k^{r+2p} $ such that
$$dP_k^{r+2p}=P_{k-1}^{r+2p+1}\mbox{\ \ and\ \ }d^*P_k^{r+2p}=0.$$
For $d^*P_k^{r+2p+1}=0$ we can take $\tilde P^{r+2p+2}_k\in\cP_k^{r+2p+2}$ such that
$$d^*\tilde P_k^{r+2p+2}=-P_{k-1}^{r+2p+1}\mbox{\ \ and\ \ }d\tilde P_k^{r+2p+2}=0.$$\\
(ii) For $dP_k^{r+2p+3}=0$ we can find $\bar P^{r+2p+2}_k\in\cP_k^{r+2p+2}$ such that
$$d\bar P_k^{r+2p+2}=P_{k-1}^{r+2p+3}\mbox{\ \ and\ \ }d^*\bar P_k^{r+2p+2}=0.$$\\
For $d^*P_k^{r+2p+3}=0$ we can take $\tilde P_k^{r+2p+4}\in\cP_k^{r+2p+4} $ such that
$$d^*\tilde P_k^{r+2p+4}=-P_{k-1}^{r+2p+3}\mbox{\ \ and\ \ }d\tilde P_k^{r+2p+4}=0.$$
Define $P^{r+2p+2}_k=\tilde P^{r+2p+2}_k+\bar P^{r+2p+2}_k.$
By (i) and (ii), obviously,
$$d^*P_k^{r+2p+2}=-P_{k-1}^{r+2p+1}\mbox{\ \ and\ \ }d P_k^{r+2p+2}=P_{k-1}^{r+2p+3}.$$
By induction, we can thus construct, for each $j=p,\ldots,q,$ a~form $P_k^{r+2j}\in\cP_k^{r+2j}$
such that
$$d^*P_k^{r+2j}=-P_{k-1}^{r+2j-1}\mbox{\ \ and\ \ }d P_k^{r+2j}=P_{k-1}^{r+2j+1}$$
where $P_{k-1}^{r+2p-1}=0$ and $P_{k-1}^{r+2q+1}=0.$ Then the form
$$P_k=P^{r+2p}_k+P^{r+2p+2}_k+\cdots+P^{r+2q}_k$$ has the required properties.
\end{proof}

%%%%%%%%%%%%%%%%%%%%%%%%%%%%%%%%%%%%%%%%%%%%%%%%%%%%%%%%%%%%%%%%%%%%%%%%%%%

Not too much is known in general about the spaces $MT^{(r,p,q)}_k.$
In the case $\bR^3,$ a~basis  for the space $MT^{(1,0,0)}_k=H^1_k$ is given in \cite{leu}
and orthonormal bases for the spaces $MT^{(1,0,0)}_k=H^1_k$ and $MT^{(0,0,1)}_k$ are constructed in \cite{cac}.
In the case $\bR^4,$
a~procedure has been worked out in \cite{del} for constructing bases for the spaces $MT^{(r,p,q)}_k.$
In the general case $\bR^m,$
bases for the space $MT^{(1,0,0)}_k=H^1_k$ have been given in \cite{del07} and \cite{zei}.
Furthermore, denote
by $\cM^+_k$ the real vector space of left monogenic polynomials in $\RR^m$ which are homogeneous of degree $k$ and which take values in the even part $\RR_{0,m}^+$ of the Clifford algebra
$\RR_{0,m}$ and put $n=[\frac m2].$
Then the following result is well-known (see \cite{BRDS} and \cite{DSS}).

\begin{lem}\label{dimM+k}
$MT^{(0,0,n)}_k\simeq\cM^+_k$ and
$\dim MT^{(0,0,n)}_k=\dim \cM^+_k=c(k,m)$
where
\begin{equation}\label{ckm}
c(k,m)=2^{m-1}{k+m-2\choose m-2}.
\end{equation}
\end{lem}

\section{The $O(m)$-modules $H^s_k$}
In this section, we shall use some known facts from representation theory of the group
$O(m).$
We refer to \cite[5.2.2]{GW} for more details.
The classification of irreducible $O(m)$-modules is closely related to the classification
of irreducible $SO(m)$-modules. Let us recall that the latter classification is
standardly given in terms of the highest weight of the module.
In the even dimensional case $m=2n,$ a~highest weight of an irreducible
 $SO(m)$-module is a vector $\la=(\la_1,\ldots,\la_n)$ of integers satisfying
 the relation $\la_1\geq\la_2\geq\ldots\geq\la_{n-1}\geq|\la_n|.$
 On the other hand, in the odd dimensional case  $m=2n+1$, the vector $\la$ satisfies the condition
 $\la_1\geq\ldots\geq\la_n\geq 0.$

The classification of $O(m)$-modules also differs in even and odd dimensions.
In the case when the dimension $m$ is odd, each irreducible representation $V$ of the group $O(m)$ remains irreducible even as a~representation over the special orthogonal group $SO(m).$ Moreover, the irreducible $O(m)$-representation
$V=V_{(\lambda,\epsilon)}$ is uniquely determined by the highest weight $\lambda$ of $V$ for $SO(m)$
and by a~number $\epsilon\in\{\pm 1\}.$

In the even dimensional case, the situation is a~bit more complicated.
In what follows, denote by $V_{\lambda}$ an irreducible $SO(m)$-module with the highest weight $\lambda.$
Assume now that $m=2n$ and $V$ is an irreducible $O(m)$-module. Then there are two possibilities. The first one is that the module $V$ remains irreducible even as $SO(m)$-module and $V$ is isomorphic to $V_{\lambda}$ for some $\lambda=(\lambda_1,\ldots,\lambda_{n-1},0).$ In this case, the irreducible $O(m)$-representation
$V=V_{(\lambda,\epsilon)}$ is again uniquely determined by the highest weight $\lambda$
and a number $\epsilon\in\{\pm 1\}.$ On the other hand, there is yet another possibility. It may happen that, as $SO(m)$-module, $V$ is reducible and $V\simeq V_{\lambda}\oplus V_{\bar\lambda}$ for some $\lambda=(\lambda_1,\ldots,\lambda_n)$ with $\lambda_n>0$ and $\bar\lambda=(\lambda_1,\ldots,\lambda_{n-1},-\lambda_n).$ In that case, we denote $V$ by $V_{(\lambda,0)}.$
See \cite[5.2.2]{GW} for details.

Now let $n$ be a~positive integer such that either $m=2n$ or $m=2n+1.$
For $1\leq s\leq n$ and $k\in\NN_0,$ let $\lambda_k^s$ be the vector in $\bR^n$ having $s$ non-zero components and given by
$$\lambda_k^s=(k+1,1,\ldots,1,0,\ldots,0).$$ Moreover, put $\lambda_0^0=(0,\ldots,0)\in\RR^n.$

The following theorem tells us that all (non-trivial) $O(m)$-modules $H^s_k$ are irreducible and mutually inequivalent. Moreover, it gives a characterization of the corresponding irreducible
$O(m)$-modules using the classification mentioned above. Using the Weyl dimensional formula, we can then compute the dimension of all these spaces.
As a~corollary, a~formula for the dimension of the space $MT^{(r,p,q)}_k$ is obtained.

\begin{thm}{\label{Hodge}} The following properties hold:
\begin{itemize}
\item[(a1)] $H^0_0\simeq V_{(\lambda^0_0,1)}$ and $H^m_0\simeq V_{(\lambda^0_0,-1)}.$\\
Moreover, for $k\geq 1,$ $H^0_k\simeq \{0\}$ and $H^m_k\simeq \{0\}.$
\item[(a2)] Let $m=2n+1,$ $1\leq s\leq n$ and let $k\in\NN_0.$\\ Then
$H^s_k\simeq V_{(\lambda^s_k,1)}$ and  $H^{m-s}_k\simeq V_{(\lambda^s_k,-1)}.$
\item[(a3)] Let $m=2n,$ $1\leq s\leq n-1$ and let $k\in\NN_0.$\\ Then
$H^s_k\simeq V_{(\lambda^s_k,1)},$ $H^{m-s}_k\simeq V_{(\lambda^s_k,-1)}$ and $H^n_k\simeq V_{(\lambda^n_k,0)}.$
\item[(a4)]
Let $1\leq s\leq m-1$ and let $k\in\NN_0.$
Putting
\begin{equation}\label{dkms}
d(k,m,s)={m-2\choose s-1}{k+m-2\choose m-2}
\frac{(2k+m)(k+m-1)}{(k+s)(k+m-s)},
\end{equation}
we have that $\dim H^s_k=d(k,m,s).$
\end{itemize}
\end{thm}

\begin{cor}{\label{dimMT}}
Let $r,p$ and $q$ be non-negative integers such that $p\leq q$ and $r+2q\leq m.$ Then
$$\dim MT^{(r,p,q)}_k  =\sum_{j=p}^q d(k,m,r+2j)+\sum_{j=p}^{q-1} d(k-1,m,r+2j+1)$$
where for $1\leq s\leq m-1,$ $d(k,m,s)$ is given by (\ref{dkms}).  Furthermore,
\begin{equation*}
d(k,m,s)=
\left\{
\begin{array}{ll}
1 & \mbox{\ \ for\ }k=0\mbox{\ and\ }s=0,m;\medskip\\
0 & \mbox{\ \ for\ }k\geq 1\mbox{\ and\ }s=0,m;\medskip\\
0 &  \mbox{\ \ for\ }k=-1.
\end{array}
\right.
\end{equation*}
\end{cor}

There is a~useful additional information describing a decomposition of
harmonic forms.

\begin{lem}{\label{harmonic_forms}}
Let $0\leq s\leq m$ and let $k\in\NN_0.$ Then we have that
$$\Ker_k^s\Delta=  H^s_k\oplus U^s_k\oplus V^s_k\oplus W^s_k$$ where
$  H^s_k$ denotes the space of solutions of the Hodge-de Rham equations and 
 $U^s_k,$ $V^s_k$ and $W^s_k$ are irreducible $O(m)$-modules
with the following properties:
\begin{itemize}
\item[(a1)] $  H^0_0\simeq V_{(\lambda^0_0,1)}$ and $  H^m_0\simeq V_{(\lambda^0_0,-1)}.$\\
Moreover, for $k\geq 1,$ $  H^0_k\simeq \{0\}$ and $  H^m_k\simeq \{0\}.$
\item[(a2)] Let $m=2n+1,$ $1\leq s\leq n$ and let $k\in\NN_0.$\\ Then
$  H^s_k\simeq V_{(\lambda^s_k,1)}$ and  $  H^{m-s}_k\simeq V_{(\lambda^s_k,-1)}.$
\item[(a3)] Let $m=2n,$ $1\leq s\leq n-1$ and let $k\in\NN_0.$\\ Then
$  H^s_k\simeq V_{(\lambda^s_k,1)},$ $  H^{m-s}_k\simeq V_{(\lambda^s_k,-1)}$ and $  H^n_k\simeq V_{(\lambda^n_k,0)}.$
\item[(a4)] For $0\leq s\leq m$ and $k\geq 0,$  we have that $\dim   H^s_k=d(k,m,s).$
\item[(b)]  For $1\leq s\leq m$ and $k\geq 1,$ we have that $U^s_k\simeq   H^{s-1}_{k-1}$ and $U^0_k\simeq\{0\}.$
%\item[(b2)] $\dim U^s_k=d(k-1,m,s-1)$ for $1\leq s\leq m$ and $k\geq 1,$
\item[(c)] For $0\leq s\leq m-1$ and $k\geq 1,$ we have that $V^s_k\simeq   H^{s+1}_{k-1}$ and $V^m_k\simeq\{0\}.$
%\item[(c2)] $\dim V^s_k=d(k-1,m,s+1)$ for $0\leq s\leq m-1$ and $k\geq 1,$
\item[(d)] For $1\leq s\leq m-1$ and $k\geq 2,$ we have that
$W^s_k\simeq   H^s_{k-2},$
$W_1^s\simeq\{0\},$\\ $W^0_k\simeq\{0\}$ and $W^m_k\simeq\{0\}.$
%\item[(d2)] $\dim W^s_k=d(k-2,m,s)$ for $1\leq s\leq m-1$ and $k\geq 2.$
\end{itemize}
\end{lem}

Proofs of Theorem \ref{Hodge} and Lemma \ref{harmonic_forms} can be
easily deduced from results proved in \cite{H}.

\section{The decomposition of the kernel of the Hodge Laplacian}

 In Lemma \ref{harmonic_forms}, the kernel of the Hodge Laplacian has been decomposed into irreducible $O(m)$-modules
$$\Ker_k^s\Delta=H^s_k\oplus U^s_k\oplus V^s_k\oplus W^s_k.$$
Now we would like to give an analytic description of the spaces $U^s_k,$ $V^s_k$ and $W^s_k.$

\begin{thm}\label{UVW} Let $0\leq s\leq m$ and let $k\in\NN_0.$ Then the following properties hold:
\begin{itemize}
\item[(a)] $U_k^s\simeq\Ker_k^s\; dd^*/\Ker^s_k\; d^*$
and $V_k^s\simeq\Ker_k^s\; d^*d/\Ker^s_k\; d,$
\item[(b)] $W_k^s\simeq\Ker_k^s\Delta\; /\left(\Ker_k^s\; dd^*\cap\Ker_k^s\; d^*d\right),$
\item[(c)] $\Ker_k^s\; dd^*\cap\Ker_k^s\; d=H^s_k\oplus U^s_k$ and
$\Ker_k^s\; d^*d\cap\Ker_k^s\; d^*=H^s_k\oplus V_k^s,$
\item[(d)] $\Ker_k^s\; dd^*\cap\Ker_k^s\; d^*d=H^s_k\oplus U^s_k\oplus V^s_k.$
\end{itemize}
\end{thm}

The Fisher decomposition \cite[p. 167]{GM} tells us that
$$\cP_k=\bigoplus_{j=0}^{[k/2]}r^{2j}\cH_{k-2j}\mbox{\ \ and thus\ \ }
\cP_k^s=\bigoplus_{j=0}^{[k/2]}r^{2j}\Ker^s_{k-2j}\Delta.$$
%Here $[x]$ is the greatest integer not greater than a~real number $x.$
From Lemma \ref{harmonic_forms}, we obtain the decomposition
\begin{equation}\label{decomP}
\cP_k^s=H^s_k\oplus
\bigoplus_{j=0}^{[k/2]}r^{2j}U^s_{k-2j}\oplus
\bigoplus_{j=0}^{[k/2]}r^{2j}V^s_{k-2j}\oplus
\bigoplus_{j=0}^{[k/2]}r^{2j}Z^s_{k-2j}
\end{equation}
where $Z^s_k=r^2 H^s_{k-2}\oplus W^s_k.$

Putting
$$X^s_{k,j}=r^{2j}Z^s_{k-2j}\cap\Ker^s_k\;d\mbox{\ \ and\ \ }Y^s_{k,j}=r^{2j}Z^s_{k-2j}\cap\Ker^s_k\;d^*,$$
\medskip
%{\it Question.} $X^s_{k,j}=r^{2j}X^s_{k-2j,0}$ and  $Y^s_{k,j}=r^{2j}Y^s_{k-2j,0}$?\medskip\\
we first prove

\begin{lem}{\label{kernels}}
The following properties hold:
\begin{itemize}
\item[(a)] $r^{2j}Z^s_{k-2j}=X^s_{k,j}\oplus Y^s_{k,j}$ and $X^s_{k,j}\simeq Y^s_{k,j}\simeq H^s_{k-2-2j};$
\item[(b)] $\Ker^s_k\; d=H^s_k\oplus
\bigoplus_{j=0}^{[k/2]}r^{2j}U^s_{k-2j}\oplus
\bigoplus_{j=0}^{[k/2]}X^s_{k,j};$
\item[(c)] $\Ker^s_k\; d^*=H^s_k\oplus
\bigoplus_{j=0}^{[k/2]}r^{2j}V^s_{k-2j}\oplus
\bigoplus_{j=0}^{[k/2]}Y^s_{k,j}.$
\end{itemize}
\end{lem}

\begin{proof}
The proof will be given by induction on the degree $s.$
For $s=0,$ the statements are obvious. Assume that the statements (a), (b) and (c) of Lemma \ref{kernels} are true for some $s=0,\ldots, m-1.$
Then we need to verify them for $s+1.$
\medskip\\
$(\alpha)$ By Lemma \ref{poincare} and by using the assumption (c), we have that
$$\Ker^{s+1}_{k-1}\;d=d(\Ker^s_k\;d^*)=d(V^s_k)\oplus\bigoplus_{j=1}^{[k/2]}d(r^{2j}V^s_{k-2j})\oplus
\bigoplus_{j=0}^{[k/2]}d(Y^s_{k,j}).$$
Using the invariance of the differential $d$ and Lemma \ref{harmonic_forms},
we get that
$H^{s+1}_{k-1}=d(V^s_k),$ while for $j=0,\ldots,[(k-1)/2]$
$$X^{s+1}_{k-1,j}=d(r^{2(j+1)}V^s_{k-2-2j})
\mbox{\ \ and\ \ }r^{2j}U^{s+1}_{k-1-2j}=d(Y^s_{k,j}).$$
Of course, $X^{s+1}_{k-1,j}\simeq H^{s+1}_{k-3-2j}.$
\medskip\\
$(\beta)$
By Lemma \ref{poincare} and by $(\alpha),$ we have that
$$d^*(\Ker^{s+1}_{k+1}\;d)=
\bigoplus_{j=0}^{[(k+1)/2]}d^*(r^{2j}U^{s+1}_{k+1-2j})\oplus
\bigoplus_{j=0}^{[(k+1)/2]}d^*(X^{s+1}_{k+1,j}).$$
Using the invariance of the codifferential $d^*$ and Lemma \ref{harmonic_forms},
we get that for $j=0,\ldots,[k/2]$
$$r^{2j}V^s_{k-2j}=d^*(X^{s+1}_{k+1,j}),\
H^s_k=d^*(U^{s+1}_{k+1})\mbox{\ \ and\ \ }
Y^s_{k,j}=d^*(r^{2(j+1)}U^{s+1}_{k-1-2j}).$$
Hence, as $d^*(\Ker^{s+1}_{k+1}\;d)=d^*(\cP^{s+1}_{k+1}),$ by virtue of (\ref{decomP}),
we may conclude that
$$d^*(r^{2j}Z^{s+1}_{k+1-2j})=d^*(X^{s+1}_{k+1,j})\mbox{\ \ and thus\ \ }
r^{2j}Z^{s+1}_{k+1-2j}\simeq d^*(X^{s+1}_{k+1,j})\oplus Y^{s+1}_{k+1,j}.$$
This shows that $r^{2j}Z^{s+1}_{k+1-2j}=X^{s+1}_{k+1,j}\oplus Y^{s+1}_{k+1,j},$ which completes the proof.
\end{proof}

\begin{proof}[Proof of Theorem \ref{UVW}] The arguments used in the steps $(\alpha)$ and $(\beta)$ of the proof of
Lemma \ref{kernels} show that
$$\Ker_k^s\; dd^*=\Ker^s_k\; d^*\oplus U^s_k\mbox{\ \ and\ \ }
\Ker_k^s\; d^*d=\Ker^s_k\; d\oplus V^s_k.$$ This easily implies all statements of Theorem \ref{UVW}.
\end{proof}

\subsection*{Acknowledgment}

R. L\'avi\v cka and V. Sou\v cek acknowledge the financial support from the grant GA 201/08/0397.
This work is also a part of the research plan MSM 0021620839, which is financed by the Ministry of Education of the Czech Republic.

%%%%%%%%%%%%%%%%%%%%%%%%%%%%%%%%%%%%%%%%%%%%%%%%%%%%%%%
%%%%%%%%%%%%%%%%%%%%%%%%%%%%%%%%%%%%%%%%%%%%%%%%%%%%%%

%%%%%%%%%%%%%%%%%%%%%%%%%%%%%%%%%%%%%%%%%%%%%%%%%%%%%%%%%%%%%%%%%%%%%%%%%%%%%%

\bigskip\bigskip\bigskip

\noindent
Richard Delanghe,\\ Clifford Research Group, Department of Mathematical Analysis,\\ Ghent University, Galglaan 2, B-9000 Gent, Belgium\\
email: \texttt{richard.delanghe@ugent.be}

\bigskip

\noindent
Roman L\'avi\v cka and Vladim\'ir Sou\v cek,\\ Mathematical Institute, Charles University,\\ Sokolovsk\'a 83, 186 75 Praha 8, Czech Republic\\
email: \texttt{lavicka@karlin.mff.cuni.cz} and \texttt{soucek@karlin.mff.cuni.cz}

%%%%%%%%%%%%%%%%%%%%%%%%%%%%%%%%%%%%%%%%%%%%%%%%%%%%%%%%%%%%%%%%%%%%%%%%%%%%%%%%
%%%%%%%%%%%%%%%%%%%%%%%%%%%%%%%%%%%%%%%%%%%%%%%%%%%%%%%%%%%%%%%%%%%%%%%%%%%%%%%%

%%%%%%%%%%%%%%%%%%%%%%%%%%%%%%%%%%%%%%%%%%%%%%%%%%%%%%%%%%%%%%%%%%%%%%%%%%%%%%%%%%%%%%%%%%%%%%%%%%%%%
%%%%%%%%%%%%%%%%%%%%%%%%%%%%%%%%%%%%%%%%%%%%%%%%%%%%%%%%%%%%%%%%%%%%%%%%%%%%%%%%%%%%%%%%%%%%%%%%%%%%%%%%%%%%%%

\end{document}